\documentclass[11pt]{article}
\usepackage{german}
\usepackage[T1]{fontenc}
\usepackage[ansinew]{inputenc}
\usepackage{latexsym}
\usepackage{textcomp}
\usepackage{color}
\usepackage{amscd}
\usepackage{amsmath}
\input{amssym.def}
\input{amssym}

     \newcommand{\fa}{\goth{a}}
     \newcommand{\fb}{\goth{b}}

     \newcommand{\Ll}{\Bbb{L}}  

     \newcommand{\Q}{\Bbb{Q}}
     
     \newcommand{\Z}{\Bbb{Z}}

    \newcommand{\ol}[1]{\overline{#1}}

    \newcommand{\ti}[1]{\tilde{#1}}
    \newcommand{\group}[1]{\langle{#1}\rangle}

    \newcommand{\ot}{\otimes}
    \newcommand{\me}{^{-1}}
    \newcommand{\mal}{^{\times}}
    
    \newcommand{\df}{\stackrel{\mathrm{def}}{=}}
    \newcommand{\mr}{\mathrm}
    \newcommand{\clo}{^{\mr{c}}}

    \newcommand{\zl}{{\Bbb{Z}_l}}
    
    \newcommand{\ql}{{\Bbb{Q}_l}}

    \newcommand{\into}{\rightarrowtail}
    \newcommand{\onto}{\twoheadrightarrow}
    \newcommand{\lto}{\longrightarrow}
    \newcommand{\da}{\downarrow}

    \newcommand{\pht}{\phantom}

    \def\daz#1{#1\da\pht{#1}}

    \newcommand{\sda}{\da_{\raisebox{-1.4mm}{$\hsp{-1}\check{}$}}}

    \newcommand{\ga}{\gamma}
    \newcommand{\Ga}{\Gamma}

    \newcommand{\la}{\lambda}
    \newcommand{\La}{\Lambda}
    \newcommand{\be}{\beta}

    \newcommand{\sn}{\par\smallskip\noindent}
    \newcommand{\mn}{\par\medskip\noindent}
    \newcommand{\bn}{\par\bigskip\noindent}
    \newcommand{\bbn}{\par\bigskip\bigskip\noindent}
    
    \newcommand{\Section}[2]{\bbn {\large #1\,. \ {\sc #2}}
                             \nopagebreak
                             \nz}

    \newcommand{\nf}[2]{\\[1.5ex]
                        \bmp{1cm}
                         (#1)
                        \emp 
                        \bmp{13.5cm}
                         \bct
                          $#2$
                         \ect
                        \emp\\[1.5ex]
         }

    \newcommand{\nz}{\\[1ex]}

    \newcommand{\hsp}[1]{\hspace*{#1mm}}

    \newcommand{\mmargin}{
     \textheight 230truemm
     \textwidth 155truemm
     \topmargin -10truemm
     \oddsidemargin 5truemm
     \evensidemargin 5truemm
     }

    \newcommand{\bmp}{\begin{minipage}}
    \newcommand{\emp}{\end{minipage}}
    \newcommand{\btb}{\begin{tabular}}
    \newcommand{\etb}{\end{tabular}}
    \newcommand{\barr}{\begin{array}}
    \newcommand{\earr}{\end{array}}
    \newcommand{\bit}{\begin{itemize}}
    \newcommand{\eit}{\end{itemize}}
    \newcommand{\ben}{\begin{enumerate}}
    \newcommand{\een}{\end{enumerate}}
    \newcommand{\bct}{\begin{center}}
    \newcommand{\ect}{\end{center}}
    \newcommand{\bfr}{\begin{flushright}}
    \newcommand{\efr}{\end{flushright}}
    \newcommand{\bea}{\begin{eqnarray*}}
    \newcommand{\eea}{\end{eqnarray*}}
    \newcommand{\bqo}{\begin{quote}}
    \newcommand{\eqo}{\end{quote}}
    \newcommand{\bdc}{\begin{description}}
    \newcommand{\edc}{\end{description}}
    \newcommand{\bdia}{\begin{CD}}
    \newcommand{\edia}{\end{CD}}

    \definecolor{light}{gray}{.3}

    \newcommand{\ind}{\mathrm{ind\,}}

    \newcommand{\res}{\mathrm{res\,}}
    
    \newcommand{\Det}{\mathrm{Det\,}}
    \newcommand{\sr}[2]{{\,\stackrel{#1}{#2}\,}}
    \newcommand{\fra}[2]{{\,\frac{#1}{#2}\,}}

    \newcommand{\theorem}{\sn
           \bdc
           \item[{\sc Theorem.}] \em }

    \newcommand{\Stop}{\edc \sn\rm} 

    \newcommand{\Lemma}[1]{\sn
           \bdc
           \item[{\sc Lemma {#1}.}] \em }
    \newcommand{\proposition}{\sn
           \bdc
           \item[{\sc Proposition.}] \em }

    \newcommand{\proof}{{\sc Proof.} \ }
    \newcommand{\remark}{\sn{\sc Remark.} \ }
    
    \newcommand{\definition}{\sn
           \bdc
           \item[{\sc Definition.}] \em }

\def\lwg{{\La_\wedge G}}
\def\ab{{\mr{ab}}}
\def\lwgab{{\La_\wedge G^\ab}}
\def\lwgs{{\La_\wedge G'}}

\def\lwgab{{\La_\wedge G^\ab}}

\def\lwcgak{{\La\clo_\wedge\Ga_k}}

\def\lwcgakm{{(\La\clo_\wedge\Ga_k)\mal}}
\def\lwcgaksm{{(\La\clo_\wedge\Ga_{k'})\mal}}

\def\qq{{{\cal{Q}}}}
\def\qwg{{\qq_\wedge G}}
\def\qwgs{{\qq_\wedge G'}}
\def\qwcgak{{\qq_\wedge\clo\Ga_k}}
\def\qwcgaks{{\qq_\wedge\clo\Ga_{k'}}}
\def\defl{{\mr{defl}}}
\def\gab{{G^\ab}}
\def\deflggab{{\defl_G^\gab}}
\def\kab{{K_\ab}}

\def\lks{{\la_{K/k'}}}
\def\gs{{G'}}

\def\ver{{\mr{ver}}}
\def\cts{{{\cal{T}}'}}
\def\Res{{\mr{Res}}}
\def\Resggs{{\Res_G^\gs}}
\def\resggs{{\res\!_G^\gs}}
\def\tkk{{t_{K/k}}}
\def\tks{{t_{K/k'}}}
\def\tka{{t_{\kab/k}}}
\newcommand{\LL}{{{\mbox{{\boldmath$\mr{L}$}}}}}
\def\LLs{{\LL'}}
\def\Lls{{\Ll'}}
\def\Llab{{\Ll^\ab}}
\def\Tr{{\mr{Tr}}}
\def\tra{{\mr{tr}\!_A}}
\def\HOM{{\mr{HOM}}}
\def\homs{{\mr{Hom}^\ast}}
\def\plc{{\psi_l\chi}}
\def\plcs{{\psi_l\chi'}}
\def\chis{{\chi'}}
\def\lkk{{L_{K/k}}}
\def\lks{{L_{K/k'}}}

\def\laks{{\la_{K/k'}}}
\def\lab{{\la_{\kab/k}}}
\def\rlg{{R_lG}}
\def\rlgs{{R_lG'}}
\def\indgsg{{\ind_\gs^G}}
\def\gbws{{\fb_\wedge'}}
\def\gaw{{\fa_\wedge}}
\def\taw{{\tau(\gaw)}}
\def\gsgg{{G'/[G,G]}}
\def\lwgsgg{{\La_\wedge(\gsgg)}}
\def\xes{{x_1'}}
\def\ys{{y'}}
\def\trabs{{\tra\gbws}}
\def\trbes{{\tra\be'}}
\def\gss{{{g'}}}

\mmargin

\begin{document}

\title{Equivariant Iwasawa theory\,: an example}
 \author{Jürgen Ritter \ $\cdot$ \ Alfred Weiss \
 \thanks{We acknowledge financial support provided by NSERC and the University of Augsburg.}
 }
\date{\pht\today}

\maketitle




 \bbn
 This note is meant to justify a remark made in the
 introduction of [6] according to which the ``main conjecture'' of equivariant
 Iwasawa theory, as formulated in [2, p.564], holds when $G=G(K/k)$ is a pro-$l$ group
 with an abelian subgroup $G'$ of index $l$.
 \mn We quickly repeat the general set-up and, in doing so, refer
 the reader to [5,\S1] for facts and notation that is taken from
 our earlier papers on Iwasawa theory. Namely, $l$ is a fixed odd
 prime number and $K/k$ a Galois extension of totally real number
 fields, with $k/\Q$ and $K/k_\infty$ finite, where $k_\infty$ is
 the cyclotomic $\zl$-extension of $k$. Throughout it will be
 assumed that Iwasawa's $\mu$-invariant $\mu(K/k)$ vanishes. We
 also fix a finite set $S$ of primes of $k$ containing all primes
 above $\infty$ and all those whose ramification index in $K/k$ is
 divisible by $l$.
 \mn In this situation it is shown in [5] that the ``main
 conjecture'' of equivariant Iwasawa theory would follow from two
 kinds of hypothetical congruences between values of Iwasawa
 $L$-functions, and one of these kinds, the so-called {\em torsion
 congruences}, has meanwhile been verified in [6]. The purpose of
 the present paper is to show that the torsion congruences already
 suffice to obtain the whole conjecture in the special case when
 $G$ is a pro-$l$ group with an abelian subgroup $G'$ of index
 $l$. Before stating the precise theorem we need to recall some
 notation (compare [5,\S1]).
 \ben\item[] $\lwg$ is the $l$-completion of the localization
 $\La_\bullet G$ which is obtained from the Iwasawa algebra
 $\La G=\zl[[G]]$ by inverting all central elements which are
 regular in $\La G/l\La G$; $\qwg$ is the total
 ring of fractions of $\lwg$;
 \item[] $T(\qwg)=\qwg/[\qwg,\qwg]$ is the quotient of
 $\qwg$ by Lie commutators;
 \item[] if $G$ is a pro-$l$ group, then (see [5,\S2] \footnote{$\rlg$ is the ring of all (virtual) ${\Q_l}\!\clo$-characters
 of $G$ with open kernel; $\Ga_k=G(k_\infty/k)\,;\ \lwcgak={\zl}\!\!\clo\ot_{\Z_l}\La_\wedge\Ga_k$ with
 ${\zl}\!\clo$ the ring of integers in a fixed algebraic closure ${\Q_l}\!\clo$ of ${\ql}$})
 \nf{LD}{\barr{ccc} K_1(\lwg)&\sr{\Ll}{\lto}&T(\qwg)\\
 \daz{\Det}&&\daz{\substack{\Tr\\ \simeq}}\\
 \HOM(\rlg,\lwcgakm)&\sr{\LL}{\lto}&\homs(\rlg,\qwcgak)\earr}
 is the logarithmic diagram defining the {\em logarithmic
 pseudomeasure}
 $$\tkk\in T(\qwg)\quad\mr{by}\quad\Tr(\tkk)=\LL(\lkk)$$
 where $\lkk=L_{K/k,S}\in\HOM(\rlg,\lwcgakm)$ is the Iwasawa $L$-function.
 \een
 \theorem With $K/k$ and $S$ as at the beginning and $G=G(K/k)$ a
 pro-$l$ group, $\tkk$ is integral (i.e., $\tkk\in T(\lwg)$)
 whenever $G$ has an abelian subgroup $G'$ of index $l$.\Stop
 As a corollary, by [5, Proposition 3.2] and [6, Theorem],
 $\lkk\in\Det K_1(\lwg)$\,, which implies the conjecture (see [3, Theorem A]),
 up to its uniqueness assertion. However, $SK_1(\qq G)=1$ because
 each simple component, after tensoring up with a suitable
 extension field of its centre, becomes isomorphic to a matrix
 ring of dimension a divisor of $l^2$ by the proof of [2,
 Proposition 6], as the character degrees $\chi(1)$ all divide
 $l$. Now apply [7, p.334, Corollary].
 \mn The proof of the theorem is carried out in \S2; before, in a
 short \S1, we introduce restriction maps  $$\Res_G^\gs:T(\qwg)\to
 T(\qwgs)\quad\mr{and}\quad\Res_G^\gs:\homs(\rlg,\qwcgak)\to\homs(\rlgs,\qwcgaks)$$
 making the diagram
 $$\barr{ccccc}K_1(\lwg)&\sr{\Ll}{\lto}&T(\qwg)&\sr{\Tr}{\lto}&\homs(\rlg,\qwcgak)\\
 \daz{\res_G^\gs}&&\daz{\Res_G^\gs}&&\daz{\Res_G^\gs}\\
 K_1(\lwgs)&\sr{\Lls}{\lto}&T(\qwgs)&\sr{\Tr'}{\lto}&\homs(\rlgs,\qwcgaks)\earr$$
 commute for any pair of pro-$l$ groups $G=G(K/k)$ and
 $G'=G(K/k')\le G$ such that $[G:G']$ is finite. We remark that replacing $\Res_G^\gs$ by the
 ``natural'' restriction map,
 $$(\res_G^\gs f)(\chis)=f(\ind_\gs^G\chis)\,,\
 f\in\homs(\rlg,\qwcgak)\,,\ \chis\in\rlgs\,,$$ does not work, because induction and Adams
 operations do not commute.
 \Section{1}{Res}
 Let $G=G(K/k)$ be a pro-$l$ group and $G'=G(K/k')\le G$ an open
 subgroup. Recall that $\Psi:\lwcgak\to\lwcgak$ is
 the map induced by $\Psi(\ga)=\ga^l$ for $\ga\in \Ga_k$ (compare
 [5,\S1]) and that $\psi_l$ is the $l^\mr{th}$ Adams operation on
 $R_l(-)$.
 \definition\quad
 $\Res_G^\gs:\homs(\rlg,\qwcgak)\to\homs(\rlgs,\qwcgaks)$ \ sends
 $f$ to $$\Res_G^\gs f=\big[\,\chis\mapsto
 f(\ind_\gs^G\chis)+\sum_{r\ge1}\fra{\Psi^r}{l^r}(f(\psi_l^{r-1}\chi))\,\big]
 \quad\footnote{A similar definition regarding
 $H_0(G;\ql\!\clo G)$ for finite $l$-groups $G$ appears in
 [1].}\,,$$
 where \ $\chi=\psi_l(\ind_\gs^G\chis)-\ind_\gs^G(\plcs)$\,. \Stop
 To justify the definition we must show that the sum
 $\sum_{r\ge1}$ is actually a finite sum. For this, let $\{t\}$ be
 a set of coset representatives of $G'$ in $G$, so $G=\dot\cup_t
 t\gs$, and define
 $$m(g)=\min\{r\ge0:g^{l^r}\in\gs\}\quad\mr{for}\quad g\in G\,.$$ Then
 $$\ind_\gs^G\chis(g)=\sum_t\dot{\chis}(g^t)=\sum_{\{t:m(g^t)=0\}}\chis(g^t)\,,$$
 if, as usual, $\dot\chis$ coincides with $\chis$ on $G'$ and vanishes
 on $G\setminus\gs$. Hence,
 $$\barr{l}\chi(g)=(\ind_\gs^G\chis)(g^l)-\ind_\gs^G(\psi_l\chis)(g)\\
 =\sum_{m(g^{lt})=0}\chis(g^{lt})-\sum_{m(g^t)=0}\chis(g^{lt})
 =\sum_{m(g^t)=1}\chis(g^{lt})\,.\earr$$
 If $r_0$ is such that $G^{l^{r_0}}\subset G'$, then
 $\psi_l^{r_0-1}\chi=0$ and $\sum_{r\ge1}=\sum_{r=1}^{r_0-2}$\,,
 because the sum $\sum_{m(g^t)=1}$ is empty when $g\in
 G^{l^{r_0-1}}$\,.
 \sn It remains to show that $\Res_G^\gs
 f\in\homs(\rlgs,\qwcgaks)$\,, i.e., $\Res_G^\gs f$ is a Galois
 stable homomorphism, compatible with W-twists (see [5,\S1]), and
 taking values in $\qwcgaks$. The first property is easily checked
 and the third follows from the second as in [2, proof of Lemma
 9]. We turn to twisting.
 \sn Let $\rho'$ be a type-W character of $G'$, so $\rho'$ is
 inflated from $\Ga_{k'}$, and write $\rho'=\res_G^\gs\rho$ with
 $\rho$ inflated from $\Ga_k$ to $G$. Then
 $$f(\ind_\gs^G(\rho'\chis))=f(\rho\cdot\ind_\gs^G\chis)=\rho^\sharp
 (f(\ind_\gs^G\chis))=(\rho')^\sharp (f(\ind_\gs^G\chis))$$
 as $f(\ind_\gs^G\chis)\in\qwcgaks$. Moreover, since $\psi_l$ is
 multiplicative,
 $$\psi_l(\ind_\gs^G(\rho'\chis))-\ind_\gs^G(\psi_l(\rho'\chis))=\psi_l(\rho\cdot\ind_\gs^G\chis)-\ind_\gs^G((\rho')^l\cdot\psi_l\chis)=\rho^l\cdot\chi$$
 and thus
 $$\barr{l}\fra{\Psi^r}{l^r}(f(\psi_l^{r-1}(\rho^l\cdot\chi))))=\fra{\Psi^r}{l^r}f(\rho^{l^r}\cdot\psi_l^{r-1}\chi)=\\
 \fra{\Psi^r}{l^r}((\rho^{l^r})^\sharp(f(\psi_l^{r-1}\chi)))=\rho^\sharp(\fra{\Psi^r}{l^r}
 f(\psi_l^{r-1}\chi))=(\rho')^\sharp(\fra{\Psi^r}{l^r}
 f(\psi_l^{r-1}\chi))\,.\earr$$
 \Lemma{1} The diagram below commutes. In it, $\LL$ and $\LLs$ are
 the lower horizontal maps of the logarithmic diagram (LD) for
 $G$ and $G'$, respectively.\nz
 {\rm(HD)}\hsp{20}$\barr{ccc}\HOM(\rlg,\lwcgakm)&\sr{\LL}{\lto}&\homs(\rlg,\qwcgak)\\
 \daz{\res_G^\gs}&&\daz{\Res_G^\gs}\\
 \HOM(\rlgs,\lwcgaksm)&\sr{\LLs}{\lto}&\homs(\rlgs,\qwcgaks)\,.\earr$
 \Stop
 \sn Indeed, for $f\in\HOM(\rlg,\lwcgakm)$ we get
 $$\barr{l}(\Res_G^\gs\LL f)(\chis)=(\LL f)(\ind_\gs^G\chis)+\sum_{r\ge1}\fra{\Psi^r}{l^r}[(\LL f)(\psi_l^{r-1}\chi)]\\
 \dot{=}(\LL f)(\ind_\gs^G\chis)+\sum_{r\ge1}\fra{\Psi^r}{l^r}[\log(f(\psi_l^{r-1}\chi))-
 \fra{\Psi}{l}\log(f(\psi_l^r\chi))]\\
 =(\LL f)(\ind_\gs^G\chis)+
 \sum_{r\ge1}\fra{\Psi^r}{l^r}\log(f(\psi_l^{r-1}\chi))-\sum_{r\ge2}\fra{\Psi^r}{l^r}\log(f(\psi_l^{r-1}\chi))\\
 =(\LL f)(\ind_\gs^G\chis)+\fra{\Psi}{l}\log(f(\chi))=
 \fra1l\log\fra{f(\ind\chis)^l}{\Psi(f(\psi_l\ind\chis))}+\fra{\Psi}{l}\log\fra{f(\psi_l\ind\chis)}{f(\ind\psi_l\chis)}\\[1.5mm]
 =\fra1l\log\fra{f(\ind\chis)^l\cdot\Psi f(\psi_l\ind\chis)}{\Psi
 f(\psi_l\ind\chis)\cdot\Psi
 f(\ind\psi_l\chis)}=\fra1l\log\fra{f(\ind\chis)^l}{\Psi
 f(\ind\psi_l\chis)}\\[1mm]
 =(\LLs\res_G^\gs f)(\chis)\ .\earr $$
 The dotted equality sign, $\dot=$, is due to the
 congruence \ $\fra{f(\chi)^l}{\Psi f(\plc)}\equiv1\mod l\lwcgak$
 (see [5,\S1]) and to $\chi(1)=0$, so $(\psi_l^{r-1}\chi)(1)=0$ for
 every $r$. In fact, with $\ti\chi\df\psi_l^{r-1}\chi$, we have  $$f(\ti\chi)^l\equiv\Psi f(\psi_l\ti\chi)\mod l\lwcgak \implies
 f(\ti\chi)^{l^s}\equiv\Psi^sf(\psi_l^s\ti\chi)\equiv\Psi^sf(\ti\chi(1)1)=1\mod
 l\lwcgak$$ for big enough $s$. Thus $(\LL
 f)(\ti\chi)=\log(f(\ti\chi))-\fra{\Psi}{l}\log(f(\psi_l\ti\chi))$
 as `log' converges on an element a power of which is $\equiv
 1\mod l\lwcgak$.
 \sn The proof of Lemma 1 is complete.
 \mn By means of the trace isomorphism $\Tr:T(-)\to\homs(-)$ we
 next transport $\Res_G^\gs$ to $\Res_G^\gs:T(\qwg)\to T(\qwgs)$,
 i.e., the diagram
 \nf{TD}{\barr{ccc}T(\qwg)&\sr{\Tr}{\lto}&\homs(\rlg,\qwcgak)\\
 \daz{\Res_G^\gs}&&\daz{\Res_G^\gs}\\
 T(\qwgs)&\sr{\Tr'}{\lto}&\homs(\rlgs,\qwcgaks)\earr}
 commutes.
 \Lemma{2}\quad $\barr{ccc}K_1(\lwg)&\sr{\Ll}{\to}&T(\qwg)\\
 \daz{\res_G^\gs}&&\daz{\Res_G^\gs}\\
 K_1(\lwgs)&\sr{\Lls}{\to}&T(\qwgs)\earr$ \quad commutes and\quad
 $\Res_G^\gs\tkk=\tks$\,.\Stop
 The first claim follows from gluing together the diagrams (LD),
 (HD), (TD) and applying [2, Lemma 9]; the second claim follows
 from $\res_G^\gs\lkk=\lks$ [2, Proposition 12].
 \mn The next lemma already concentrates on the case when $G'$ is
 abelian and $[G:G']=l$. We set $A=G/G'=\group{a}$ and observe
 that $a$ acts on $G'$ by conjugation.
 \Lemma{3} Let $\tau:\lwg\to T(\lwg)$ denote the canonical map and
 $g\in G$. If $G'$ is abelian \footnote{whence $\tau':\lwgs\to
 T(\lwgs)$ is the identity map}
 and of index $l$ in $G$, then
 $$\Resggs(\tau g)=\left\{\barr{ll}\sum_{i=0}^{l-1}g^{a^i}& \mr{if}\ g\in\gs\\
 g^l&\mr{if}\ g\notin\gs\,.\earr\right.$$\Stop
 To see this, we apply $\Tr'$ to both sides and employ the formula
 \ $\Tr'(\tau'g)(\chis)=\chis(g)\ol g$ with $\ol g$ denoting the
 image of $g\in\gs$ in $\Ga_{k'}$ (see [5,\S1])\,:
 \ben\item $(\Tr'\Resggs(\tau g))(\chis)=\Resggs(\Tr(\tau
 g))(\chis)=\Tr(\tau g)(\ind_\gs^G\chis)+\fra{\Psi}{l}\Tr(\tau
 g)(\chi)$ since $G^l\subset G'$. Now, if $g\in G'$,
 $\Tr(\tau g)(\indgsg\chis)=\sum_{i=0}^{l-1}\chis(g^{a^i})\ol
 g$ and $\chi(g)=0$. On the other hand, if $g\notin\gs$, $\Tr(\tau
 g)(\indgsg\chis)=0$ and $\fra{\Psi}{l}\Tr(\tau
 g)(\chi)=\fra1l\indgsg\chis(g^l)\ol{g}^l=\chis(g^l)\ol{g}^l$
 since we may choose $a=g\mod\gs$.
 \item
 $\Tr'(\sum_{i=0}^{l-1}g^{a^i})(\chis)=\sum_{i=0}^{l-1}\chis(g^{a^i})\ol
 g$, since $g^{a^i}$ and $g$ have the same image in $\Ga_k$ and so
 in $\Ga_{k'}$. On the other hand,
 $\Tr'(g^l)(\chis)=\chis(g^l)\ol{g^l}$\,.
 \een
 The lemma is established.
 \remark The lemma has two immediate generalizations. Firstly, if
 $\Ga\,(\simeq\zl)$ is a central subgroup of $G$ contained in
 $G'$, then the elements of $T(\lwg)$ can uniquely be written as
 $\sum_g\be_g\tau(g)$ with $\be_g\in\La_\wedge\Ga$ and $g$ running
 through a set of preimages of conjugacy classes of $G/\Ga$ (see
 [3, Lemma 5]). For each summand we have
 $$\Resggs(\be_g\tau(g))=\left\{\barr{ll}\sum_{i=0}^{l-1}\be_gg^{a^i}& \mr{if}\ g\in\gs\\
 \Psi(\be_g)g^l&\mr{if}\ g\notin\gs\,.\earr\right.$$
 Secondly, if $G'$ is no longer abelian, then the equality in the
 lemma has to be replaced by
 $$\Resggs(\tau g)=\left\{\barr{ll}\sum_{i=0}^{l-1}\tau'(g)^{a^i}& \mr{if}\ g\in\gs\\
 \tau'(g^l)&\mr{if}\ g\notin\gs\,.\earr\right.$$
 \Section{2}{Proof of the theorem}
 In this section $G=G(K/k)$ is a pro-$l$ group and $G'=G(K/k')$ an
 abelian subgroup of index $l$ ($K/k$ is as in the introduction).
 As before, $A=G/\gs=\group{a}$\,, and we set $\hat
 A=1+a+\cdots+a^{l-1}$\,.  \mn If $G$ itself is abelian, the theorem holds by [4,\S5,
 Example 1], whence we assume that $G$ is non-abelian.
 \Lemma{4} Assume that there exists an element $x\in T(\lwg)$ such
 that $\deflggab x=\deflggab\tkk$ and $\Resggs x=\Resggs\tkk$\,.
 Then $\tkk\in T(\lwg)$.\Stop
 Denoting by $\kab$ the fixed field of the finite group $[G,G]$,
 we first observe, because of [5, Lemma 2.1] and Lemma 2, that
 $\deflggab\tkk=\tka$ and $\Resggs\tkk=\tks$ are integral\,:
 indeed, a logarithmic pseudomeasure is integral whenever the
 group is abelian.
 \sn From [4, Proposition 9] we obtain a power $l^n$ of $l$ such that
 $l^n\tkk\in T(\lwg)$. Consider the element $\ti x=l^n(x-\tkk)\in
 T(\lwg)$. It satisfies $\deflggab\ti x=0=\Resggs\ti x$\,. We are
 going to prove $\ti x=0$ which implies $x=\tkk$ because
 $\homs(\rlg,\qwcgak)$, and so $T(\qwg)$, is torsionfree; whence
 the the lemma will be verified.
 \sn The proof of $\ti x=0$ employs the commutative diagram shown
 in the proof of [5, Proposition 2.2]\,:
 $$\barr{cccccc}1+\gaw&\into&(\lwg)\mal&\sr{\deflggab}{\onto}&(\lwgab)\mal&\\
 \sda&&\daz{\Ll}&&\daz{\Llab}&\\
 \tau(\gaw)&\into& T(\lwg)&\sr{\deflggab}{\onto}&\lwgab&,\earr$$
 in which $\Ll$ is extended to $(\lwg)\mal$ by means of the
 canonical surjection $(\lwg)\mal\onto K_1(\lwg)$ and
 $\gaw=\ker(\lwg\to\lwgab)$. The diagram yields
 a $v\in(\lwg)\mal$ with $\Ll(v)=\ti x$, simply because
 $\deflggab\ti x=0$. Combining diagram (HD) of Lemma 1 and that
 appearing in Lemma 2, we arrive at $$\LLs(\resggs(\Det
 v))=\Resggs(\LL(\Det v))=\Resggs(\Tr\,\Ll(v))=\Tr'(\Resggs\ti
 x)=0$$ and, with $\resggs$ replaced by $\deflggab$, at
 $$\LL^\ab(\deflggab(\Det v))=\deflggab(\LL(\Det
 v))=\deflggab(\Tr\,\Ll(v))=\Tr^\ab(\deflggab\ti x)=0\,,$$ since
 $\LL$ and $\Tr$ commute with deflation.
 \sn The first displayed formula in [3, p.46] now implies that
 $\resggs(\Det v)$ and $\deflggab(\Det v)$ are torsion elements in
 $\HOM(\rlgs,\lwcgaksm)$ and $\HOM(R_l(\gab),\lwcgakm)$,
 respectively. Moreover, the first paragraph of the proof of [5,
 Proposition 3.2] therefore shows that $\Det v$ itself is a
 torsion element in $\HOM(\rlg,\lwcgakm)$. Consequently, for some
 natural number $m$, $(\Det v)^{l^m}=1$, so $l^m\LL(\Det
 v)=0=l^m\Tr(\Ll v)=\Tr(l^m\ti x)$\,, and $\ti x=0$ follows, as has
 been claimed.
 \bbn We now introduce the commutative diagram
 $$\barr{ccccc}\taw&\into& T(\lwg)&\onto&\lwgab=T(\lwgab)\\
 \daz{\Res}&&\daz{\Resggs}&&\daz{\Res}\\
 \gbws=\tau'(\gbws)&\into&\lwgs=T(\lwgs)&\onto&\lwgsgg=T(\lwgsgg)\earr$$
 with exact rows (of which the upper one has already appeared in
 the diagram shown in the proof of the preceding lemma). The images
 of all vertical maps are fixed elementwise by $A$ because of
 Lemma 3. Thus we can turn the diagram into
 \nf{D}{\barr{ccccc}\taw&\into&T(\lwg)&\onto&\lwgab\\
 \daz{\Res}&&\daz{\Resggs}&&\daz{\Res}\\
 \gbws^A&\into&(\lwgs)^A&\to&(\lwgsgg)^A\\
 \sda&&\sda&&\sda\\ \hat H^0(A,\gbws)&\to&\hat H^0(A,\lwgs)&\to&
 \hat H^0(A,\lwgsgg)\earr}
 with exact rows and canonical lower vertical maps.
 \Lemma{5} In (D), the left vertical column is exact and the left
 bottom horizontal map is injective.\Stop
 \proof The ideal $\gaw$ is (additively) generated by the elements
 $g(c-1)$ with $g\in G,\,c\in[G,G]$\,; those with $g\in\gs$ generate
 $\gbws$. We compute $\Resggs\tau(g(c-1))$, using Lemma 3\,:
 \ben\item if $g\in G'$,
 $\Resggs\tau(g(c-1))=\sum_{i=0}^{l-1}((gc)^{a^i}-g^{a^i})=\sum_{i=0}^{l-1}\Big(g(c-1)\Big)^{a^i}\in\trabs$\,,
 \item if $g\notin\gs$,
 $\Resggs\tau(g(c-1))=\Resggs(\tau(gc)-\tau(g))=(gc)^l-g^l=g^lc^{\hat
 A}-g^l=0$\,, since
 \nf{$\star$}{[G,G]^{\hat A}=1} by $[G,G]\dot=[G,G']$ and $[G,G']^{\hat A}=((G')^{a-1})^{\hat A}=1$ as $(a-1)\hat
 A=0$. Here, the dotted equality sign, $\dot=$, results from the
 equation
 $$\barr{l}[bg_1',b^ig_2']=(g_1')\me b\me (g_2')\me b^{-i}bg_1' b^ig_2'=(g_1')\me
 ({g_2'}\me)^{b}(g_1')^{b^i}g_2'=\\ \Big((g_1')\me(g_1')^{b^i}\Big)\Big(({g_2'}\me)^{b}g_2'\Big)
 \in[G',G]\cdot[G,G']\le[G,G']\earr$$ for $g_1',g_2'\in G'$ and $b\in G\setminus\gs$,  because $G'$ is abelian and normal in
 $G$.\een
 Thus, $\Resggs\taw=\trabs$\,, which proves the first claim of the
 lemma.
 \sn The second claim follows from $\hat
 H\me(A,\lwgsgg)=0$ and this in turn from the trivial action of
 $A$ on $\gsgg$ and the torsion freeness of $\lwgsgg$.
 \sn Lemma 5 is established.
 \bn As seen in diagram (D), there is an element $x_1\in
 T(\lwg)$ with $\deflggab x_1=\tka$. We define $\xes\in\lwgs$ by
 $\Resggs x_1=\tks+\xes$. Because of [5, Lemma 3.1], $\xes$ is
 fixed by $A$. We want to change $x_1$ modulo $\taw$ so that
 the new $\xes$ becomes zero: then we have arrived at an $x\in
 T(\lwg)$ as assumed in Lemma 4 and the theorem will have been
 confirmed.
 \sn The above change is possible if, and only if,
 $\xes\in\Resggs(\taw)$ and so, because of Lemma 5, if $\xes$ is in
 $\cts\df\tra(\lwgs)$\,, the $A$-trace ideal of the $A$-action on
 $\lwgs$.
 \proposition $\xes\in\cts$ is achievable.\Stop
 This is seen as follows.  From [5,\S1] we recall the existence of
 {\em pseudomeasures} $\lab,\laks$ in $K_1(\lwgab)$ and
 $K_1(\lwgs)$, respectively, satisfying
 $\Det\lab=L_{\kab/k},\,\Det\laks=\lks$ (so
 $\Llab(\lab)=\tka,\,,\Lls(\laks)=\tks$). From [5, 2.~of Proposition
 3.2] and [6, Theorem] we know that
 $$\fra{\ver(\lab)}{\laks}\equiv1\mod\cts$$
 where `ver' is the map induced from the transfer homomorphism
 $\gab\to\gs$.
 \sn Let $y\in(\lwg)\mal$ have $\deflggab y=\lab$ and set $\resggs
 y=\laks\cdot\ys$. Then $$\ys=\fra{\resggs
 y}{\laks}\equiv\fra{\ver(\lab)}{\laks}\equiv1\mod\cts$$
 (see the proof of [5, Proposition 3.2]). Moreover,
 $\ys\in1+\gbws$. Now, $x_1\df\Ll(y)$ has
 $\Resggs x_1=\Resggs\Ll(y)=\tks+\xes$ with $\xes\df\Lls(\ys)$, and $\xes\in\gbws$
 because of the commutativity of
 $$\barr{cccccc}1+\gbws&\into&(\lwgs)\mal&\onto&(\lwgsgg)\mal&\\
 \daz{\Lls}&&\daz{\Lls}&&\daz{\ol{\Lls}}&\\ \gbws&\into&\lwgs&\onto&\lwgsgg&.\earr$$
 Hence, the proposition (and therefore the theorem) will be proved, if
 \nf{2.0}{\xes=\Lls(y')\in\cts.}
 However, Lemma 5 gives
 $$\ys\in(1+\gbws^A)\cap(1+\cts)=1+(\gbws^A\cap\cts)=1+\trabs$$
 and as $\Lls(\ys)=\fra1l\log\fra{\ys^l}{\Psi(\ys)}$ (compare
 [3, p.39]),
 we see that \nf{2.1}{\Lls(\ys)\in\cts\quad\mr{if}\quad
 \fra{\ys^l}{\Psi(\ys)}\equiv1\mod l\cts\,.}
 So it suffices to show this last congruence.
 \sn Write $y'=1+\trbes$ with $\be'\in\gbws$. Since
 $(1+\trbes)^l\equiv1+(\trbes)^l\mod l\cts$, the congruence in (2.1) is
 equivalent to \nf{2.2}{(\trbes)^l\equiv\Psi(\trbes)\mod l\cts.}
 On picking a central subgroup $\Ga\,(\simeq\zl)$ of $G$ and
 writing \ $\be'=\sum_{g',c}\be_{g',c}\,g'(c-1)$ with elements
 $\be_{g',c}\in\La_\wedge\Ga,\,g'\in G',\,c\in[G,G]$, we obtain
 \nf{2.2a}{\barr{l}(\trbes)^l=\Big(\sum_{g',c}\be_{g',c}\tra(g'(c-1))\Big)^l
 \equiv\sum_{g',c}(\be_{g',c})^l\Big(\tra(g'(c-1))\Big)^l\\ \equiv\sum_{g',c}\Psi(\be_{g',c})\Big((\tra(g'c))^l-(\tra
 g')^l\Big)\mod l\cts\earr} and
 \nf{2.2b}{\Psi(\trbes)=\sum_{g',c}\Psi(\be_{g',c})\Big(\tra((g'c)^l)-\tra({g'}^l)\Big)}
 as $\Psi$ and $\tra$ commute. Thus congruence (2.2) will
 result from Lemma 6 below,
 since then subtracting (2.2b) from (2.2a)
 yields the sum
 $$\barr{l}\sum_{g',c}\Psi(\be_{g',c})\Big((\tra(g'c))^l-\tra((g'c)^l)-(\tra
 g')^l+\tra({g'}^l)\Big)\equiv\\
 \sum_{g',c}\Psi(\be_{g',c})\Big(-l(g'c)^{\hat
 A}+l{g'}^{\hat
 A}\Big)\equiv\sum_{g',c}(-l)\Psi(\be_{g',c}){g'}^{\hat
 A}(c^{\hat A}-1)\equiv0\mod l\cts\,,\earr$$
 by $(\star)$ of the proof of Lemma 5.
 \Lemma{6} \quad $(\tra g')^l-\tra(
 {g'}^l)\equiv -l{g'}^{\hat A}\mod l\cts$ \ for $\gss\in\gs$\,.\Stop
 \proof Set \ $\ti A=\Z/l\times A$ \ and make \
 $M=\mr{Maps}(\Z/l,A)$ \
 into an $\ti A$-set by defining $m^{(z,a^i)}(x)=m(x-z)\cdot a^i$. Then
 $$(\tra\gss)^l=(\sum_{i=0}^{l-1}\gss^{a^i})^l=\sum_{m\in
 M}\prod_{z\in\Z/l}\gss^{m(z)}=\sum_{m\in M}\gss^{\sum_{z\in\Z/l}m(z)}$$ with $\sum_zm(z)$ read in
 $\Z[A]$.
 \sn We compute the subsums of $\sum_m$ in which $m$ is
 constrained to an $\ti A$-orbit.
 \sn If $m\in M$ has stabilizer $\{(0,1)\}$ in $\ti A$, then the
 $\ti A$-orbit sum is $$\barr{l}\sum_{(z,a^i)\in\ti
 A}\gss^{\sum_{v\in\Z/l}m^{(z,a^i)}(v)}=\sum_{(z,a^i)\in\ti
 A}\gss^{\sum_{v\in\Z/l}m(v-z)a^{i}}=\\ \sum_{(z,a^i)}\gss^{\sum_vm(v)a^i}=l\sum_i(\gss^{\sum_vm(v)})^{a^i}=l\cdot\tra(\gss^{\sum_vm(v)})\in
 l\cts\ .\earr$$
 Note that no $m\in M$ is stabilized by $(0,a^i)$ with
 $a^i\neq1$\,: for $m(z)=m^{(0,a^i)}(z)=m(z)a^i$ implies $a^i=1$.
 It follows that the stabilizers of the elements with stabilizer
 different from $\{(0,1)\}$ must be cyclic of order $l$ and
 different from $\{(0,a^i):0\le i\le l-1\}$ and therefore $=\group{(1,a^j)}$ for
 a unique $j\mod l$.
 \sn One now checks that for each $j$ there is exactly one $\ti
 A$-orbit with stabilizer $\group{(1,a^j)}$ and that it
 is represented by $m_j\,,\,m_j(z)=a^{jz}$. Moreover,
 $\{(0,a^i):0\le i\le l-1\}$ is a transversal of the stabilizer of
 $m_j$ in $\ti A$.
 \sn For each $j$, the sum of $\gss^{\sum_zm(z)}$ over the $\ti
 A$-orbit of $m_j$ is
 $\sum_i\gss^{\sum_zm_j^{(0,a^i)}(z)}=\sum_i\gss^{\sum_za^{jz}a^i}$\,.
 If $j=0$, this is
 $\sum_i\gss^{la^i}=\tra(\gss^l)$\,,
 accounting for that term in the claim. If $j\neq0$, it is
 $\sum_i\gss^{\hat A\cdot a^i}=l\gss^{\hat A}$, and summing over
 $j\neq0$ gives $(l-1)l\cdot\gss^{\hat A}\equiv-l\cdot\gss^{\hat A}\mod l\cts$
 because $l^2\cdot\gss^{\hat A}=l\cdot\tra(\gss^{\hat A})\in l\cts$.
 \sn This finishes the proof of Lemma 6.
 \bbn
 {\sc References}
 \bn
 \btb{rp{13cm}}
 \,[1]    & Oliver, R.~and Taylor, L.R., {\em Logarithmic
           descriptions of Whitehead groups and class groups for
           $p$-groups.} Memoirs of the AMS (1988), vol.~76, no.~392\\
 \,[2-6] & Ritter, J.~and Weiss, A.\,, \newline
           \hsp{2}2. {\em Towards equivariant Iwasawa theory, II.}
           Indag.\,Mathemat.~{\bf 15} (2004), 549-\hsp{6}572\newline
           \hsp{2}3. {\em $\cdots$, III.} Math.\,Ann.~{\bf 336}
           (2006), 27-49\newline
           \hsp{2}4. {\em $\cdots$, VI.} Homology, Homotopy and
           Applications {\bf 7} (2005), 155-171\newline
           \hsp{2}5. {\em Non-abelian pseudomeasures and congruences
           between abelian Iwasawa $L$-\hsp{6}functions.} To appear in
           Pure and Applied Math Quarterly (2007)\newline
           \hsp{2}6. {\em Congruences between abelian pseudomeasures.}
           Preprint (2007)\\
 \,[7]   & Wang, S., {\em On the commutator group of a simple algebra.}
           Amer.\,J.\,Math.~{\bf 72} (1950), 323-334

 \etb
 \nopagebreak[4]
 \mn {\footnotesize \bct Institut für Mathematik $\cdot$
 Universität Augsburg $\cdot$ 86135 Augsburg $\cdot$ Germany \\
 Department of Mathematics $\cdot$ University of Alberta $\cdot$
 Edmonton, AB $\cdot$ Canada T6G 2G1   \ect

\end{document}